\documentclass{article}

\newcommand{\be}{\begin{equation}}
\newcommand{\ee}{\end{equation}}
\newcommand{\bea}{\begin{eqnarray}}
\newcommand{\eea}{\end{eqnarray}}
\newcommand{\barray}{\begin{array}}
\newcommand{\earray}{\end{array}}
\newcommand{\pa}{\partial}
\newcommand{\nn}{\nonumber}
\newcommand{\bitem}{\begin{itemize}}
\newcommand{\eitem}{\end{itemize}}
\newtheorem{teo}{Theorem}[section]
\newcommand{\bt}{\begin{teo}}
\newcommand{\et}{\end{teo}}
\newtheorem{Def}{Definition}[section]
\newcommand{\bd}{\begin{Def}}
\newcommand{\ed}{\end{Def}}
\newtheorem{lem}{Lemma}[section]
\newcommand{\bl}{\begin{lem}}
\newcommand{\el}{\end{lem}}
\newtheorem{prop}{Proposition}[section]
\newcommand{\bp}{\begin{prop}}
\newcommand{\ep}{\end{prop}}
\newtheorem{cor}{Corollary}[section]
\newcommand{\bc}{\begin{cor}}
\newcommand{\ec}{\end{cor}}
\newtheorem{ex}{Example}[section]
\newcommand{\bex}{\begin{ex}}
\newcommand{\eex}{\end{ex}}
\newtheorem{rem}{Remark}[section]
\newcommand{\br}{\begin{rem}}
\newcommand{\er}{\end{rem}}

\catcode `\@=11
\@addtoreset{equation}{section}

\begin{document}

\begin{center}
{\Large \textbf{Compatible metrics of constant \\ Riemannian
curvature: local geometry, \\ nonlinear equations
and integrability\footnote{This work was supported by
the Alexander von Humboldt Foundation (Germany),
the Russian Foundation for Basic Research
(grant No. 99--01--00010) and the INTAS
(grant No. 99--1782).}}}
\end{center}

\bigskip
\bigskip

\centerline{\large {O. I. Mokhov}}
\bigskip
\bigskip

\section{Introduction} \label{vved}

In the present paper, the nonlinear equations
describing all the nonsingular pencils of
metrics of constant Riemannian curvature are derived and
the integrability of these nonlinear equations
by the method of inverse scattering problem
is proved. These results were announced
in our previous paper \cite{1}.
For the flat pencils of metrics the corresponding
statements and proofs were presented in the present author's
work \cite{2}, \cite{3}, where the method of integrating
the nonlinear equations for the nonsingular flat pencils of
metrics was proposed. In \cite{4} the Lax pair for
the nonsingular flat pencils of metrics was demonstrated.
This Lax pair is generalized to the case of
arbitrary nonsingular pencils of metrics of
constant Riemannian curvature (see examples and
interesting applications in \cite{4}).

In this paper, it is proved that
all the nonsingular pairs of compatible
metrics of constant Riemannian curvature are described
by special integrable reductions of nonlinear equations
defining orthogonal curvilinear coordinate systems
in the spaces of constant curvature.
Note that the problem of description for
the pencils of metrics of constant Riemannian curvature
is equivalent to the problem of description
for compatible nonlocal Poisson brackets of hydrodynamic type
generated by metrics of constant Riemannian curvature
(compatible Mokhov--Ferapontov brackets \cite{5})
playing an important role in the theory
of systems of hydrodynamic type.

Recall that two pseudo-Riemannian contravariant metrics
$g_1^{ij} (u)$ and $g_2^{ij} (u)$ are called {\it compatible}
if for any linear combination of these metrics
$g^{ij} (u) = \lambda_1 g_1^{ij} (u) + \lambda_2 g_2^{ij} (u)$,
where $\lambda_1$ and $\lambda_2$ are arbitrary
constants for which $\det ( g^{ij} (u) ) \not\equiv 0$,
the coefficients of the corresponding
Levi-Civita connections and the components of
the corresponding tensors of Riemannian curvature are related
by the same linear formula:
$\Gamma^{ij}_k (u) = \lambda_1 \Gamma^{ij}_{1, k} (u) +
\lambda_2 \Gamma^{ij}_{2, k} (u)$ and
$R^{ij}_{kl} (u) = \lambda_1 R^{ij}_{1, kl} (u)
+ \lambda_2 R^{ij}_{2, kl} (u)$
(in this case, we shall say also that {\it the metrics
$g^{ij}_1 (u)$ and $g^{ij}_2 (u)$ form a pencil of metrics}) \cite{1}.
Flat pencils of metrics, that is nothing but
 compatible nondegenerate local Poisson brackets
of hydrodynamic type (compatible Dubrovin--Novikov
brackets \cite{6}), were introduced in \cite{7}.
Two pseudo-Riemannian contravariant metrics
$g_1^{ij} (u)$ and $g_2^{ij} (u)$ of constant Riemannian
curvature
$K_1$ and $K_2$ respectively are called {\it compatible}
if any linear combination of these metrics
$g^{ij} (u) = \lambda_1 g_1^{ij} (u) + \lambda_2 g_2^{ij} (u),$
where $\lambda_1$ and $\lambda_2$ are
arbitrary constants for which
$\det ( g^{ij} (u) ) \not\equiv 0$, is
a metric of constant Riemannian curvature
$\lambda_1 K_1 + \lambda_2 K_2$
and the coefficients of the corresponding
Levi-Civita connections are related by the same
linear formula:
$\Gamma^{ij}_k (u) = \lambda_1 \Gamma^{ij}_{1, k} (u) +
\lambda_2 \Gamma^{ij}_{2, k} (u)$ \cite{1}.
In this case, we shall also say that
{\it the metrics
$g_1^{ij} (u)$ and $g_2^{ij} (u)$ form a pencil of
metrics of constant Riemannian curvature} \cite{1}.
It is obvious that all these definitions are mutually
consistent, so that if compatible metrics are
metrics of constant Riemannian curvature, then
they form a pencil of metrics of constant Riemannian curvature,
and if compatible metrics are flat, then they form
a flat pencil of metrics. A pair of pseudo-Riemannian
metrics $g_1^{ij} (u)$ and $g_2^{ij} (u)$ is called
{\it nonsingular} if the eigenvalues of this pair
of metrics, that is, the roots of the equation
$\det ( g_1^{ij} (u) -  \lambda g_2^{ij} (u)) =0,$
are distinct (a pencil of metrics which is
formed by a nonsingular pair of metrics is also
called {\it nonsingular}). In \cite{1}, it is proved
that an arbitrary nonsingular pair of metrics is
compatible if and only if there exist local coordinates
$u = (u^1,...,u^N)$ such that both the metrics
are diagonal in these coordinates and have the
following special form (one can consider that
one of the metrics, here
$g^{ij}_2 (u),$ is an arbitrary
diagonal metric):
$g^{ij}_2 (u) = g^i (u) \delta^{ij}$ and
$g^{ij}_1 (u) = f^i (u^i) g^i (u) \delta^{ij},$
where $f^i (u^i),$ $i=1,...,N,$ are functions of
single variable (generally speaking, complex).
It is obvious that the eigenvalues of the
considered pair of metrics are given by the functions
$f^i (u^i),$ $i=1,...,N.$

\section{Equations for nonsingular
pairs of compatible \\  metrics of
constant Riemannian curvature}    \label{p2}

Consider the problem of description
for nonsingular pairs of compatible metrics of
constant Riemannian curvature.
As is proved in \cite{1} (see also Introduction above),
it is sufficient to classify pairs of diagonal
metrics of constant Riemannian curvature of
the special type
$g_2^{ij} (u) = g^i (u) \delta^{ij}$ and
$g_1^{ij} (u) = f^i (u^i) g^i (u) \delta^{ij},$
where $f^i (u^i),$ $i= 1,...,N,$ are
arbitrary (generally speaking, complex)
functions of single variable,
$g^i (u),$ $i=1,...,N,$ are arbitrary functions.

The problem of description for $N \times N$
diagonal metrics of constant Riemannian curvature $K$, that is,
metrics $g_2^{ij} (u) = g^i (u) \delta^{ij}$  such that
their tensor of Riemannian curvature
has the form $R^{ij}_{kl} =
K (\delta^i_l \delta^j_k - \delta^i_k \delta^j_l),$ where
$K = {\rm const},$ is the classical problem of
differential geometry.
This problem is equivalent to the problem
of description of orthogonal curvilinear coordinate
systems in $N$-dimensional spaces of constant curvature $K$.
Recently Zakharov showed that the Lam\'e equations
describing orthogonal curvilinear coordinate systems in
Euclidean or pseudo-Euclidean $N$-dimensional
spaces are integrated by the method
of inverse scattering problem \cite{8}
(see also an algebraic-geometric approach in \cite{9}).
Similarly, the equations describing
orthogonal curvilinear coordinate systems in
$N$-dimensional spaces of constant curvature
are integrated by the method of inverse scattering problem.

Introduce the standard classical notation
\bea
&&
g^i (u) = {\varepsilon^i \over (H_i (u))^2 },\ \ \ d\, s^2 =
\sum_{i=1}^N  \varepsilon^i
(H_i (u))^2 (d u^i)^2, \\
&&
\beta_{ik} (u) = {1 \over H_i (u)} {\pa H_k \over \pa u^i},\ \ \
i \neq k, \label{vra}
\eea
where $H_i (u)$ are {\it the Lam\'e coefficients}
and
$\beta_{ik} (u)$ are {\it the rotation coefficients},
$\varepsilon^i = \pm 1,$ $i = 1,...,N$.
Although, in our case, all the functions
are, generally speaking, complex, we shall use
formulae, which are convenient for using also in
the purely real case.

\bt \label{syst}
Nonsingular pairs of compatible
metrics of constant Riemannian
curvature $K_1$ and $K_2$ are described
by the following consistent integrable nonlinear systems:
\be
{\pa \beta_{ij} \over \pa u^k}
=\beta_{ik} \beta_{kj},\ \ \ i\neq j,\ \ i\neq k,\ \ j\neq k, \label{lam1}
\ee
\be
\varepsilon^i {\pa \beta_{ij} \over \pa u^i}+
\varepsilon^j {\pa \beta_{ji} \over \pa u^j}+\sum_{s\neq i,\
s\neq j} \varepsilon^s \beta_{si} \beta_{sj} =
- K_2 H_i H_j,\ \ \ i\neq j, \label{lam2}
\ee
\bea
&&
\varepsilon^i f^i (u^i) {\pa  \beta_{ij}
 \over \pa u^i}+  {1 \over 2} \varepsilon^i (f^i)' \beta_{ij} +
\varepsilon^j f^j (u^j)
{\pa
\beta_{ji} \over \pa u^j}+
{1 \over 2} \varepsilon^j (f^j)' \beta_{ji} + \nn\\
&&
\sum_{s\neq i,\
s\neq j} \varepsilon^s f^s (u^s) \beta_{si} \beta_{sj} =
- K_1 H_i  H_j,\ \ \ i\neq j, \label{lam3}
\eea
where $f^i (u^i),$ $i=1,...,N,$ are arbitrary given
functions of one variable.
\et

\br
Equations (\ref{lam1}) are the well-known
$n$-wave equations (the Darboux equations \cite{10}),
and equations (\ref{lam1}), (\ref{lam2}) describe
orthogonal curvilinear coordinate systems in
$N$-dimensional spaces of constant curvature
$K_2,$ in particular, for $K_2 = 0$ we get the famous
Lam\'e equations. Equations (\ref{lam3})
define a nontrivial nonlinear differential reduction
of equations (\ref{lam1}), (\ref{lam2}).
\er

\br
If there are coinciding nonzero constants
among the functions $f^i (u^i)$, then equations
(\ref{lam1})--(\ref{lam3}) also describe a pencil of
metrics of constant Riemannian curvature, but, obviously, this
pencil is not nonsingular.
\er

The system of equations (\ref{vra})--(\ref{lam3})
is consistent for any functions $f^i (u^i)$ what can be easily
checked by direct calculations.
Moreover, the general solution of the nonlinear
system (\ref{vra})--(\ref{lam3}) depend on $N^2$
arbitrary functions of one variable (the functions
$\beta_{ij} (u)$ and $H_j (u)$
can be arbitrarily defined
on the $j$th coordinate line ).

Let us consider the conditions that the diagonal metrics
$g^{ij}_2 (u) = g^i (u) \delta^{ij}$ and
$g^{ij}_1 (u) = f^i (u^i) g^i (u) \delta^{ij},$
where $f^i (u^i),$ $i=1,...,N,$ are arbitrary
functions of one indicated variable (these
functions must not be only equal to zero identically),
are metrics of constant Riemannian curvature
$K_2$ and $K_1$ respectively.

Recall that for any diagonal metric
$\Gamma^i_{jk} (u) =0$ if all the indices $i, j, k$
are distinct. It is also obvious that
$R^{ij}_{kl} (u) = 0$ if all the indices $i, j, k, l $
are distinct. In addition, by virtue of the well-known
symmetries of the tensor of Riemannian curvature,
we have:
$$R^{ii}_{kl} (u) = R^{ij}_{kk} (u) =0,
\ \ \ \ R^{ij}_{il} (u) = -R^{ij}_{li} (u) = R^{ji}_{li} (u) =
- R^{ji}_{il} (u).$$

Thus, it is sufficient to consider the condition
$R^{ij}_{kl} (u) = K_2 (\delta^i_l \delta^j_k -
\delta^i_k \delta^j_l)$ (the necessary and sufficient
condition that the corresponding metric is a metric
of constant Riemannian curvature $K_2$) only
for the following components of the tensor
of Riemannian curvature:
 $R^{ij}_{il} (u)$,
where $i \neq j,$ $\ i \neq  l$:
$$R^{ij}_{il} = - K_2 \delta^j_l, \ \ \ \ i \neq j,\ \ i\neq l.$$

For an arbitrary diagonal metric
$g^{ij}_2 (u) = g^i (u) \delta^{ij}$, we get
$$\Gamma^i_{2, ik} (u) = \Gamma^i_{2, ki} (u) =
- {1 \over 2 g^i (u)} {\pa g^i \over \pa u^k}, \ \ \
{\rm \ for \  all\ } i, k; $$
$$\Gamma^i_{2, jj} (u) = {1 \over 2} {g^i (u) \over
(g^j (u))^2 } {\pa g^j \over \pa u^i},\ \ \ i \neq j,$$

\bea
&&
R^{ij}_{2, il} (u) = g^i (u) R^j_{2, iil} (u) = \nn\\
&&
g^i (u) \left ( {\pa \Gamma^j_{2, il} \over \pa u^i}  -
{\pa \Gamma^j_{2, ii}  \over \pa u^l} + \sum_{s=1}^N
\Gamma^j_{2, si} (u) \Gamma^s_{2, il} (u) -
\sum_{s=1}^N \Gamma^j_{2, sl} (u) \Gamma^s_{2, ii} (u) \right ).
\eea

It is necessary to consider two cases separately.

1) $j \neq l$.

$$R^{ij}_{il} =0, \ \ \ i\neq j,\ \ i\neq l,\ \ j\neq l.$$

\bea
&&
R^{ij}_{2, il} (u) =
g^i (u) \left (
- {\pa \Gamma^j_{2, ii}  \over \pa u^l} +
\Gamma^j_{2, ii} (u) \Gamma^i_{2, il} (u) -
\Gamma^j_{2, jl} (u) \Gamma^j_{2, ii} (u) -
\Gamma^j_{2, ll} (u) \Gamma^l_{2, ii} (u) \right )
=  \nn\\
&&
- {1 \over 2} g^i (u)  {\pa \over \pa u^l}
\left ( {g^j (u) \over (g^i (u))^2 } {\pa g^i \over \pa u^j }\right )
 - {1 \over 4} {g^j (u) \over (g^i (u))^2 }
{\pa g^i \over \pa u^j} {\pa g^i \over \pa u^l} \nn\\
&&
+ {1 \over 4 g^i (u)}
{\pa g^i \over \pa u^j} {\pa g^j \over \pa u^l}
- {1 \over 4} {g^j (u) \over g^i (u) g^l (u) }
{\pa g^l \over \pa u^j} {\pa g^i \over \pa u^l}= 0 .\label{r1}
\eea

Equations (\ref{r1}) are equivalent to the equations
\be
 {\pa^2 H_i \over \pa u^j \pa u^k} =
{1 \over H_j (u)} {\pa H_i \over \pa u^j} {\pa H_j \over \pa u^k}
+ {1 \over H_k (u)} {\pa H_k \over \pa u^j} {\pa H_i \over \pa u^k},
\ \ \ i\neq j,\ \ i\neq k,\ \ j\neq k,
\ee
which are equivalent, in turn, to equations (\ref{lam1}).

2) $j=l$.

$$R^{ij}_{ij} = - K_2, \ \ \ i\neq j.$$

\bea
&&
R^{ij}_{2, ij} (u) =
g^i (u) \left (
 {\pa \Gamma^j_{2, ij}  \over \pa u^i} -
{\pa \Gamma^j_{2, ii}  \over \pa u^j} +
\Gamma^j_{2, ii} (u) \Gamma^i_{2, ij} (u) + \right. \nn\\
&&
\left. \Gamma^j_{2, ji} (u) \Gamma^j_{2, ij} (u)
 -  \sum_{s=1}^N
\Gamma^j_{2, sj} (u) \Gamma^s_{2, ii} (u) \right )
=  \nn\\
&&
- {1 \over 2} g^i (u)  {\pa \over \pa u^i}
\left ( {1 \over g^j (u) } {\pa g^j \over \pa u^i }\right )
- {1 \over 2} g^i (u)  {\pa \over \pa u^j}
\left ( {g^j (u) \over (g^i (u))^2 } {\pa g^i \over \pa u^j }\right )
- {1 \over 4} {g^j (u) \over (g^i (u))^2 }
{\pa g^i \over \pa u^j} {\pa g^i \over \pa u^j} \nn\\
&&
+ {1 \over 4} {g^i (u) \over (g^j(u))^2}
{\pa g^j \over \pa u^i} {\pa g^j \over \pa u^i}
- {1 \over 4 g^j (u)}
{\pa g^j \over \pa u^i} {\pa g^i \over \pa u^i} +
\sum_{s \neq  i} {1 \over 4} {g^s (u) \over g^i (u) g^j (u) }
{\pa g^j \over \pa u^s} {\pa g^i \over \pa u^s}= - K_2.\label{r2}
\eea

Equations (\ref{r2}) are equivalent to the equations
\bea
&&
\varepsilon^i {\pa \over \pa u^i} \left ( {1 \over H_i (u)}
{\pa H_j  \over \pa u^i} \right ) +
\varepsilon^j {\pa \over \pa u^j} \left ( {1 \over H_j (u)}
{\pa H_i  \over \pa u^j} \right ) + \nn\\
&&
\sum_{s \neq i,\ s \neq j}
{\varepsilon^s \over (H_s (u))^2}
{\pa H_i \over \pa u^s} {\pa H_j \over \pa u^s}= - K_2 H_i H_j,
\ \ \ i \neq j,
\eea
which are equivalent, in turn, to equations (\ref{lam2}).

The condition that the metric
$g_1^{ij} (u) = f^i (u^i) g^i (u) \delta^{ij}$
is also a metric of constant Riemannian curvature
gives $N(N-1)/2$ additional equations (\ref{lam3}).
Note that in this case components (\ref{r1})
of the tensor of Riemannian curvature vanish automatically,
and the condition (\ref{r2}) gives the corresponding
$N(N-1)/2$ equations.

Actually, for the metric
$g^{ij}_1 (u) = f^i (u^i) g^i (u) \delta^{ij}$,
the Lam\'e coefficients and the rotation coefficients
have the form
\bea
&&
\widetilde H_i (u) = {H_i (u)  \over \sqrt {\epsilon^i f^i (u^i)}},
\ \ \ f^i (u^i) g^i (u) = {\epsilon^i \varepsilon^i \over
(\widetilde H_i (u))^2},\ \ \ \epsilon^i = \pm 1, \label{redu1}\\
&&
\widetilde \beta_{ik} (u) =
{1 \over \widetilde H_i (u)} {\pa \widetilde H_k \over
\pa u^i} = \nn\\
&&
{\sqrt {\epsilon^i f^i (u^i)} \over
\sqrt {\epsilon^k f^k (u^k)}} \left ( {1 \over H_i (u)} {\pa H_k \over
\pa u^i}  \right ) =  {\sqrt {\epsilon^i f^i (u^i)} \over
\sqrt {\epsilon^k f^k (u^k)}}
\beta_{ik} (u), \  i \neq k. \label{reduction}
\eea
Accordingly, equations (\ref{lam1}) are
automatically satisfied also for the rotation coefficients
$\widetilde \beta_{ik} (u)$, and equations
(\ref{lam2}) for $\widetilde \beta_{ik} (u)$ give equations
(\ref{lam3}).

Note that, for our further purposes, it is more convenient
to write the system
(\ref{lam1})--(\ref{lam3}) namely in this form
in order to
emphasize the reduction (\ref{redu1}), (\ref{reduction})
playing an important role in our method of integrating
this system (see the next section),
although it is easy to show that equations
(\ref{lam2}), (\ref{lam3}) for nonsingular pairs of metrics
(that is, all the functions $f^i (u^i)$ must be distinct
also in the case if they are constants) are equivalent
to the following equations (in particular,
it is more convenient to use
these equations for checking the consistency of system
(\ref{vra})--(\ref{lam3}))
\bea
&&
{\pa \beta_{ij} \over \pa u^i}
= {1 \over 2} {(f^i (u^i))' \over (f^j (u^j) - f^i (u^i))}
\beta_{ij} + {\varepsilon^i \varepsilon^j \over 2}
{(f^j (u^j))' \over (f^j (u^j) - f^i (u^i))}
\beta_{ji} - \nn\\
&&
- \sum_{s \neq i, s \neq j} \varepsilon^i \varepsilon^s
{(f^j (u^j) - f^s (u^s)) \over (f^j (u^j) - f^i (u^i))}
\beta_{si} \beta_{sj} + \varepsilon^i
{K_1 - K_2 f^j (u^j) \over (f^j (u^j) - f^i (u^i))} H_i  H_j,
 \ \ i \neq j.
\eea

\bex
If all the functions $f^i (u^i)$ are arbitrary
distinct nonzero constants:
$f^i (u^i) = c^i,$ $c^i = {\rm const} \neq 0,$
$i=1,...,N,$ $c^i \neq c^j$ for $i \neq j,$
then the integrable system
(\ref{vra})--(\ref{lam3}) takes the following form:
\bea
&&
{\pa \beta_{ij} \over \pa u^k}
=\beta_{ik} \beta_{kj},\ \ \ i\neq j,\ \ i\neq k,\ \ j\neq k,
\label{lam1co} \\
&&
{\pa \beta_{ij} \over \pa u^i}
= - \sum_{s \neq i, s \neq j} \varepsilon^i \varepsilon^s
{(c^j - c^s) \over (c^j  - c^i)}
\beta_{si} \beta_{sj} + \varepsilon^i
{K_1 - K_2 c^j  \over (c^j - c^i)} H_i  H_j,
 \ \ i \neq j, \label{lam2co} \\
&&
{\pa H_j \over \pa u^i} =
\beta_{ij} H_i,\ \ \
i \neq j. \label{vra1}
\eea
\eex

In the flat case, for $K_1=K_2=0,$ system (\ref{lam1co}),
(\ref{lam2co}) was considered in \cite{4}
(it is related to a curious triple of pairwise commuting
Monge--Amp\`ere equations). So the system
of equations (\ref{lam1co})--(\ref{vra1})
is also of special interest for applications.

\section{Compatible flat metrics and
the Zakharov method \\ of differential reductions} \label{p3}

In this section, we demonstrate the method of
integrating the system of equations
(\ref{lam1})--(\ref{lam3}) in the flat case (for $K_1 =K_2=0$)
(see \cite{2}, \cite{3}), that is, the following system (in this
section $\varepsilon^i = 1,$ $i=1,...,N$ in all formulae):
\be
{\pa \beta_{ij} \over \pa u^k}
=\beta_{ik} \beta_{kj},\ \ \ i\neq j,\ \ i\neq k,\ \ j\neq k, \label{lam1a}
\ee
\be
{\pa \beta_{ij} \over \pa u^i}+
 {\pa \beta_{ji} \over \pa u^j}+\sum_{s\neq i,\
s\neq j}  \beta_{si} \beta_{sj} =
0,\ \ \ i\neq j, \label{lam2a}
\ee
\bea
&&
 f^i (u^i) {\pa  \beta_{ij}
 \over \pa u^i}+  {1 \over 2} (f^i)' \beta_{ij} +
 f^j (u^j)
{\pa
\beta_{ji} \over \pa u^j}+
{1 \over 2}  (f^j)' \beta_{ji} + \nn\\
&&
\sum_{s\neq i,\
s\neq j}  f^s (u^s) \beta_{si} \beta_{sj} =
0,\ \ \ i\neq j, \label{lam3a}
\eea
where $f^i (u^i),$ $i=1,...,N,$ are arbitrary given
functions of one variable.

Recall the Zakharov method for integrating
the Lam\'{e} equations
(\ref{lam1a}) and (\ref{lam2a}) \cite{8}.

We must choose a matrix function $F_{ij} (s, s', u)$
with some special properties and
solve the linear integral equation
\be
K_{ij}(s, s', u)  = F_{ij}(s, s', u) + \int_s^{\infty}
\sum_l K_{il} (s, q, u) F_{lj} (q, s', u) dq.  \label{int}
\ee
Then we obtain a one-parameter family of solutions of
the Lam\'{e} equations by the
formula
\be
\beta_{ij} (s, u) = K_{ji} (s, s, u).  \label{resh}
\ee

In particular,
if $F_{ij} (s, s', u) = f_{ij} (s-u^i, s'-u^j),$
where $f_{ij} (x, y)$ is an arbitrary matrix function of two variables,
then formula (\ref{resh}) produces solutions of equations (\ref{lam1a}).
To construct solutions of
equations (\ref{lam2a}), Zakharov proposed to impose
on the ``dressing matrix function'' $F_{ij} (s-u^i, s' -u^j)$ a certain
additional linear differential relation. If
the function $F_{ij} (s - u^i, s' - u^j)$
satisfies the Zakharov differential relation
(we shall present it below), then
the rotation coefficients $\beta_{ij} (u)$ additionally satisfy
equations
(\ref{lam2a}).

Let us describe our scheme of integrating all the system
of equations
(\ref{lam1a})--(\ref{lam3a}).

\bl
If both the function $F_{ij} (s-u^i, s' -u^j)$ and
the function

\be
\widetilde F_{ij} (s-u^i, s' -u^j) =
{\sqrt {f^j (u^j -s')} \over \sqrt {f^i (u^i -s)} }
F_{ij} (s-u^i, s' -u^j)   \label{f}
\ee
satisfy the Zakharov differential relation, then
the corresponding rotation coefficients $\beta_{ij} (u)$
(\ref{resh}) satisfy all the equations (\ref{lam1a})--(\ref{lam3a}).
\el
 Actually, if $K_{ij} (s, s', u)$ is the
corresponding to the
function $F_{ij} (s-u^i, s' -u^j)$
 solution of
the linear integral equation (\ref{int}), then
\be
\widetilde K_{ij} (s, s', u) =
{\sqrt {f^j (u^j -s')} \over \sqrt {f^i (u^i -s)} }
K_{ij} (s, s', u)
\ee
is the corresponding to
function (\ref{f}) solution of equation (\ref{int}).
It is easy to prove multiplying the integral equation
(\ref{int}) by
$\sqrt {f^j (u^j -s')}/\sqrt {f^i (u^i -s)}: $

\be
\widetilde K_{ij}(s, s', u)  = \widetilde F_{ij}(s -u^i, s' - u^j)
+ \int_s^{\infty}
\sum_l \widetilde K_{il} (s, q, u)
\widetilde F_{lj} (q - u^l,  s'- u^j) dq.  \label{int2}
\ee
Then both $\widetilde \beta _{ij} (s, u) = \widetilde K_{ji} (s, s, u)$
and $\beta_{ij} (s, u) = K_{ji} (s, s, u)$
satisfy the Lam\'{e} equations (\ref{lam1a}) and (\ref{lam2a}).
Besides, we
have
\be
\widetilde \beta _{ij} (s, u) = \widetilde K_{ji} (s, s, u)
=
{\sqrt {f^i (u^i -s)} \over \sqrt {f^j (u^j -s)} }
K_{ji} (s, s, u) =
{\sqrt {f^i (u^i -s)} \over \sqrt {f^j (u^j -s)} }
\beta_{ij} (s, u).
\ee

Thus, in this case, the rotation coefficients $\beta_{ij} (u)$
satisfy exactly all equations (\ref{lam1a})--(\ref{lam3a}), that is,
they generate the corresponding compatible flat metrics.

The Zakharov reduction is given by the following
differential relations \cite{8}:

\be
{\pa F_{ij} (s, s', u) \over \pa s'} +
{\pa F_{ji} (s', s, u) \over \pa s} = 0. \label{za1}
\ee

To resolve these relations for the matrix
function $F_{ij} (s- u^i, s' -u^j)$,
we can introduce $N(N-1)/2$ arbitrary functions of two variables
$\Phi_{ij} (x, y),$ $i < j$, and put for $i<j$
\bea
&&
F_{ij} (s - u^i, s' - u^j) =
{\pa \Phi_{ij} (s-u^i, s' - u^j) \over \pa s},\nn\\
&&
F_{ji} (s - u^i, s' - u^j) =
- {\pa \Phi_{ij} (s' -u^i, s - u^j) \over \pa s}, \label{zare1}
\eea
and, besides, put for any $i$
\be
F_{ii} (s - u^i, s' - u^i) = {\pa \Phi_{ii} (s-u^i, s' - u^i)
\over \pa s}, \label{zare2}
\ee
where $\Phi_{ii} (x, y)$, $i=1,...,N,$
are arbitrary skew-symmetric functions:
\be
\Phi_{ii} (x, y) = - \Phi_{ii} (y, x),
\ee
see \cite{8}.

For the function
\be
\widetilde F_{ij} (s-u^i, s' -u^j) =
{\sqrt {f^j (u^j -s')} \over \sqrt {f^i (u^i -s)} }
F_{ij} (s-u^i, s' -u^j)
\ee
the Zakharov differential relations (\ref{za1})
give exactly $N(N-1)/2$ linear partial
differential equations of the second order
 for
$N(N-1)/2$  functions of two variables
$\Phi_{ij} (s- u^i, s' - u^j), \ i < j$:
\bea
&&
{\pa \over \pa s'} \left (
{\sqrt {f^j (u^j -s')} \over \sqrt {f^i (u^i -s)} }
{\pa \Phi_{ij} (s-u^i, s' - u^j) \over \pa s}
\right ) -\nn\\
&&
{\pa \over \pa s} \left (
{\sqrt {f^i (u^i -s)} \over \sqrt {f^j (u^j -s')} }
{\pa \Phi_{ij} (s-u^i, s' - u^j) \over \pa s'} \right ) = 0,
\ \ \ \ i<j,
\eea
which are equivalent to the equations
\bea
&&
2 {\pa^2 \Phi_{ij} (s - u^i, s' - u^j) \over
\pa u^i \pa u^j} \left (f^i (u^i -s) - f^j (u^j - s') \right ) +
{\partial \Phi_{ij} (s- u^i, s' -u^j) \over \pa u^j}
{d f^i (u^i-s) \over d u^i}  - \nn\\
&&
{\partial \Phi_{ij} (s- u^i, s' -u^j) \over \pa u^i}
{d f^j (u^j-s') \over d u^j} =0, \ \ \ i<j.\label{ss}
\eea

It is interesting that all these equations (\ref{ss})
for the functions $\Phi_{ij} (s- u^i, s' - u^j)$
are of the same type and
coincide with the single equation
in the two-component case ($N = 2$), which was derived directly
from the conditions of vanishing the corresponding tensors
of Riemannian curvature (see \cite{1}).

Besides, for $N$ functions
$\Phi_{ii} (s-u^i, s' - u^i),$ we have also $N$ linear partial
differential
equations of the second order,
which are derived from the Zakharov differential relations
(\ref{za1}):
\be
{\pa \over \pa s'} \left (
{\sqrt {f^i (u^i -s')} \over \sqrt {f^i (u^i -s)} }
{\pa \Phi_{ii} (s-u^i, s' - u^i) \over \pa s}
\right ) +
{\pa \over \pa s} \left (
{\sqrt {f^i (u^i -s)} \over \sqrt {f^i (u^i -s')} }
{\pa \Phi_{ii} (s' -u^i, s - u^i) \over \pa s'} \right ) =0
\ee
or equivalently
\bea
&&
2 {\pa^2 \Phi_{ii} (s - u^i, s' - u^i) \over
\pa s \pa s'} \left (f^i (u^i -s) - f^i (u^i - s') \right ) - \nn\\
&&
{\partial \Phi_{ii} (s- u^i, s' -u^i) \over \pa s}
{d f^i (u^i-s') \over d s'}  +
{\partial \Phi_{ii} (s- u^i, s' -u^i) \over \pa s'}
{d f^i (u^i -s) \over d s} =0.  \label{ss2}
\eea

Any solution of the linear partial differential equations
of the second order
(\ref{ss}) and (\ref{ss2}) generates a one-pa\-ra\-me\-ter
family of solutions of system (\ref{lam1a})--(\ref{lam3a})
by linear relations and explicit formulae
(\ref{zare1}), (\ref{zare2}), (\ref{int}), (\ref{resh}).
Thus system (\ref{lam1a})--(\ref{lam3a}) is linearized.
The method of differential reductions is applicable
also to integrating arbitrary nonsingular pencils
of metrics of constant Riemannian curvature
(\ref{vra})--(\ref{lam3}).
In the next section, we present the Lax pair
with a spectral parameter for system
(\ref{vra})--(\ref{lam3}), what also gives the
possibility to integrate the system.

Note that equations (\ref{ss}), (\ref{ss2}) can be
easily integrated explicitly in the special functions.
If the functions $f^i (u^i)$ and $f^j (u^j)$
are not constants, then, as is proved in \cite{1},
for the description of all the corresponding pencils of
metrics, it is sufficient to consider the case
$f^i (u^i) = u^i,$
$f^j (u^j) = u^j,$ since in this case
the functions
$f^i (u^i)$ and $f^j (u^j)$ can be
chosen as the corresponding new local coordinates.
Accordingly equations (\ref{ss})
belong to one of the following three types:
\bea
&&
{\pa^2 F \over \pa u^1 \pa u^2} =
{1 \over 2 (u^1 - u^2)} {\pa F \over \pa u^1}
- {1 \over 2 (u^1 - u^2)} {\pa F \over \pa u^2}, \label{f1}\\
&&
{\pa^2 F \over \pa u^1 \pa u^2} =
- {1 \over 2 u^2} {\pa F \over \pa u^1}, \label{f2}\\
&&
{\pa^2 F \over \pa u^1 \pa u^2} =  0.   \label{f3}
\eea
The general solutions of equations (\ref{f2}) and (\ref{f3})
are the functions $$F (u^1, u^2) = g(u^1)/\sqrt{u^2}
+ h (u^2)\  {\rm \ \ \ and \ \ \ }\  F (u^1, u^2) = g(u^1) + h(u^2)$$
respectively, where $g (u^1)$ and $h (u^2)$ are
arbitrary functions.

After the change of variables
$u^1 = t + r,$ $u^2 = t - r,$ equation (\ref{f1})
becomes the equation
\be
{\pa^2 F \over \pa t^2} =
{\pa^2 F \over \pa r^2} +
{1 \over r}  {\pa F \over \pa r}, \label{f4}
\ee
which is well known in the classical geometry
(one more well-known Darboux equation)
and in the classical mathematical physics
(the equation of axially symmetric
oscillations of gas ($r$ is the radial coordinate)
or the equation for solutions of the wave equation
which correspond to processes of radiation on the
plane $(x, y)$, that is, solutions
depending only on the radial coordinate
$r = \sqrt{x^2 + y^2}$ and the time $t$).

In particular,
the functions
\be
v (t, r) = {\rm const} \int_{-1}^{+1} {\psi (t + x r)
\over \sqrt{1 - x^2}} dx,
\ee
where $\psi (t)$ is an arbitrary function,
satisfy to
the Darboux equation (\ref{f4})
with the initial conditions
\be
F (t, 0) = \psi (t), \ \ \  {\pa F \over
\pa r} (t, 0) = 0.
\ee
For the function $\psi (t, z) = \psi (t),$
given on the plane $(t, z),$ the function
$v (t, r)$ gives the mean value over
the circle of radius $r$ with the centre at the point
$(t, z)$.

The general solutions of equations
(\ref{f4}) can be constructed by the method
of separation of variables: $F (t, r) = S(t) R(r),$
\be
{S'' \over S} = {R'' \over R} + {1 \over r} {R' \over R},
\ee
\be
S''  = a S, \ \ \
r R'' + R' - a r R = 0, \ \ \
a = {\rm const}.
\ee
It is easy to write the general solutions in terms
of the special Bessel and modified Bessel functions or degenerate
hypergeometric functions.

\section{Lax pair for nonsingular pencils of metrics\\
of constant Riemannian curvature}

The Lax pair for nonsingular flat pencils of
metrics
(\ref{lam1a})--(\ref{lam3a}) was demonstrated by Ferapontov
in \cite{4}:
\bea
&&
{\pa \varphi_i \over \pa u^j} = \sqrt {(\lambda + f^i)
\over (\lambda + f^j)} \beta_{ij} \varphi_j, \ \ \ i \neq j,
\label{lax1} \\
&&
{\pa \varphi_i \over \pa u^i} = -
\sum_{k \neq i} \sqrt { (\lambda + f^k)
\over (\lambda + f^i)} \beta_{ki} \varphi_k,
\label{lax2}
\eea
where $\lambda$ is the spectral parameter.

Note that the linear problem for the Lam\'e
equations
(\ref{lam1a}), (\ref{lam2a}) was well known Darboux yet \cite{10},
see also, for example,
\cite{8}, \cite{9}:
\bea
&&
{\pa \varphi_i \over \pa u^j} =  \beta_{ij} \varphi_j, \ \ \ i \neq j,
\label{lax1a} \\
&&
{\pa \varphi_i \over \pa u^i} = -
\sum_{k \neq i}  \beta_{ki} \varphi_k.
\label{lax2a}
\eea
The condition of consistency for the linear system
(\ref{lax1a}), (\ref{lax2a}) defines the Lam\'e equations
(\ref{lam1a}), (\ref{lam2a}).
The Lax pair (\ref{lax1}), (\ref{lax2}) can be
easily derived from the classical linear problem
(\ref{lax1a}), (\ref{lax2a}) for the Lam\'e equations.
Actually, the system of equations
(\ref{lam1a})--(\ref{lam3a})
defining nonsingular flat pencils of metrics
is derived from the condition that the diagonal metric
$(\lambda + f^i (u^i)) g^i (u) \delta^{ij}$
is flat for any $\lambda$. It is obvious that
in this case (see (\ref{reduction}) and \cite{2},
\cite{3}) it is necessary to
change
$\beta_{ij}(u)$ to
$$\hat{\beta}_{ij} (u) = \sqrt {(\lambda + f^i (u^i)) \over
(\lambda + f^j (u^j))} \beta_{ij} (u)$$ in the Lam\'e
equations (\ref{lam1a}), (\ref{lam2a}).
Then the linear problem
(\ref{lax1a}), (\ref{lax2a}) becomes the Lax pair
(\ref{lax1}), (\ref{lax2}).

The Lax pair (\ref{lax1}), (\ref{lax2})
is generalized also to the case of arbitrary
nonsingular pencils of metrics of constant Riemannian curvature
(\ref{vra})--(\ref{lam3}) (see
examples in \cite{4}).
The Lax pair for the system (\ref{vra})--(\ref{lam3})
can be also easily derived from the linear
problem for the system (\ref{vra})--(\ref{lam2})
describing all the orthogonal curvilinear coordinate
systems in $N$-dimensional spaces of constant curvature $K_2$:
\bea
&&
{\pa \varphi_i \over \pa u^j} = \sqrt {\varepsilon^i
\varepsilon^j} \beta_{ij} \varphi_j, \ \ i \neq j, \label{l1}\\
&&
{\pa \varphi_i \over \pa u^i} = -
\sum_{k \neq i} \sqrt {\varepsilon^k \varepsilon^i}
\beta_{ki} \varphi_k
+ \sqrt{\varepsilon^i K_2}
H_i \psi, \label{l2}\\
&&
{\pa \psi \over \pa u^i} = -
 \sqrt{\varepsilon^i K_2}
H_i \varphi_i
\label{l3}
\eea
(the condition of consistency for the linear
system (\ref{l1})--(\ref{l3}) gives the equations
(\ref{vra})--(\ref{lam2})).
Actually, the system of equations (\ref{vra})--(\ref{lam3})
defining the pencils of metrics of constant Riemannian
curvature is equivalent to the condition that
the diagonal metric
$(\lambda + f^i (u^i)) g^i (u) \delta^{ij}$ is
a metric of constant Riemannian curvature $\lambda K_2 + K_1$
for any $\lambda$.
In this case (see (\ref{redu1}), (\ref{reduction}) and \cite{2},
\cite{3}), in the equations (\ref{lam1}), (\ref{lam2}),
it is necessary replace
$\varepsilon^i$ by $\varepsilon^i \epsilon^i$,
$H_i (u)$ by
$$\hat H_i (u) =
{H_i (u) \over \sqrt{ \epsilon^i (\lambda +
f^i (u^i))}},$$
 $\beta_{ij}(u)$ by
$$\hat{\beta}_{ij} (u) = \sqrt { \epsilon^i
(\lambda + f^i (u^i)) \over
 \epsilon^j (\lambda + f^j (u^j))} \beta_{ij} (u),$$
$K_2$ by $\hat K = \lambda K_2 + K_1,$
$\epsilon^i = \pm 1$.
Then the linear problem
(\ref{l1})--(\ref{l3}) becomes the Lax pair
with the spectral parameter for the system (\ref{vra})--(\ref{lam3}):
\bea
&&
{\pa \varphi_i \over \pa u^j} = \sqrt {\varepsilon^i (\lambda + f^i)
\over \varepsilon^j (\lambda + f^j)} \beta_{ij} \varphi_j, \ \ \ i \neq j,
\label{lax1b} \\
&&
{\pa \varphi_i \over \pa u^i} = -
\sum_{k \neq i} \sqrt {\varepsilon^k (\lambda + f^k)
\over \varepsilon^i (\lambda + f^i)} \beta_{ki} \varphi_k
+ \sqrt{ \lambda K_2 + K_1 \over \varepsilon^i (\lambda + f^i)}
H_i \psi,
\label{lax2b} \\
&&
{\pa \psi \over \pa u^i} = -
 \sqrt{  \lambda K_2 + K_1 \over \varepsilon^i (\lambda + f^i)}
H_i \varphi_i,
\label{lax3b}
\eea
where $\lambda$ is the spectral parameter.
The condition of consistency for the linear system
(\ref{lax1b})--(\ref{lax3b})
is equivalent to the equations (\ref{vra})--(\ref{lam3}).

\medskip

\begin{flushleft}
Centre for Nonlinear Studies,\\
L.D.Landau Institute for Theoretical Physics, \\
Russian Academy of Sciences\\
e-mail: mokhov@mi.ras.ru; mokhov@landau.ac.ru\\
\end{flushleft}

\end{document}